\documentclass[12pt]{amsart}

\usepackage{amssymb}
\usepackage{amsmath}
\evensidemargin\oddsidemargin

\begin{document}

\setcounter{page}{1}

\newtheorem{REMS}{Remarks\!\!}
\newtheorem{LEM}{Lemma\!\!}
\newtheorem{THE}{Theorem\!\!}

\renewcommand{\theTHE}{}
\renewcommand{\theLEM}{}
\renewcommand{\theREMS}{}

\newcommand{\eqnsection}{
\renewcommand{\theequation}{\thesection.\arabic{equation}}
    \makeatletter
    \csname  @addtoreset\endcsname{equation}{section}
    \makeatother}
\eqnsection

\def\a{\alpha}
\def\cL{{\mathcal{L}}} 
\def\EE{{\mathbb{E}}} 
\def\Ga{\Gamma}
\def\hk{{\kappa}}
\def\lbd{\lambda}
\def\lcr{\left[}
\def\lpa{\left(}
\def\pb{{\mathbb{P}}}
\def\rl{{\mathbb{R}}}
\def\rpa{\right)}
\def\rcr{\right]}
\def\Sa{S_\a}
\def\Un{{\bf 1}}

\def\claw{\stackrel{d}{\longrightarrow}}
\def\elaw{\stackrel{d}{=}}
\def\qed{\hfill$\square$}
                  
\title[Self-decomposability of the Fr\'echet distribution]
      {On the self-decomposability of the Fr\'echet distribution}

\author[Pierre Bosch and Thomas Simon]{Pierre Bosch and Thomas Simon}

\address{Laboratoire Paul Painlev\'e, Universit\'e Lille 1, Cit\'e Scientifique, F-59655 Villeneuve d'Ascq Cedex. {\em Emails} : {\tt pierre.bosch@ed.univ-lille1.fr, simon@math.univ-lille1.fr}}



\begin{abstract} Let $\{\Ga_t, \, t\ge 0\}$ be the Gamma subordinator. Using a moment identification due to Bertoin-Yor (2002), we observe that for every $t > 0$ and $\a\in (0,1)$ the random variable $\Ga_t^{-\a}$ is distributed as the exponential functional of some spectrally negative L\'evy process. This entails that all size-biased samplings of Fr\'echet distributions are self-decomposable and that the extreme value distribution $F_\xi$ is infinitely divisible if and only if $\xi\not\in (0,1),$ solving problems raised by Steutel (1973) and Bondesson (1992). We also review different analytical and probabilistic interpretations of the infinite divisibility of $\Ga_t^{-\a}$ for $t,\a > 0.$ 
\end{abstract}

\maketitle
 
\section{Introduction}

The extreme value theorem - see e.g. Theorem 8.13.1 in \cite{BGT} - states that non-degenerate distribution functions arising as limits of properly renormalized running maxima of i.i.d. random variables  belong to one of the families
$$F_0(x)\, =\, e^{-e^{-x}}, \; x\in\rl,\qquad\mbox{or}\qquad F_\xi (x)\, =\,\left\{\begin{array}{ll} 1- e^{-x^{1/\xi}} & \mbox{if $\xi > 0$}\\e^{-x^{1/\xi}} &  \mbox{if $\xi < 0$}\end{array} \right., \; x >0.$$
The distribution $F_0$ is known as the Gumbel distribution, whereas $F_\xi$ is called a Weibull distribution for $\xi > 0$ and a Fr\'echet distribution for $\xi < 0.$ In the following, we denote by $X_\xi$ the random variable with distribution function $F_\xi.$ Observe that
$$\frac{1- X_\xi}{\xi}\;\claw\; X_0\qquad \mbox{as $\xi\to 0,$}$$
so that the above parametrization is continuous in $\xi.$ In the present paper we are interested in the  self-decomposability (SD) of $X_\xi,$ referring e.g. to Section 15 in \cite{S1} for an account on self-decomposability. The Gumbel distribution is SD because of the identities
$$X_0\; \elaw\; -\log L \; \elaw -\a\log L \, +\, \a \log \Sa$$ 
for every $\a \in (0,1),$ where here and throughout $L$ stands for the standard exponential variable and $\Sa$ for the standard positive $\a$-stable variable - see e.g. Exercise 29.16 in \cite{S1} for a proof of the second identity. If $\xi \in (0,1)$ then the variable $X_\xi$ is not infinitely divisible (ID) and hence not SD, because of its superexponential distribution tails - see e.g. Theorem 26.1 in \cite{S1}. When $\xi \ge 1$, the variable $X_\xi$ has a completely monotone density and is ID by Goldie's criterion - see e.g. Theorem 4.2 in \cite{St}, or by the ME property which makes it the first-passage time of some continuous time Markov chain - see e.g. Chapter 9 in \cite{B} for an account. When $\vert\xi\vert\ge 1,$ the identity in law
$$X_\xi\; \elaw\; L^\xi$$
and the HCM theory of Thorin and Bondesson \cite{B} show that the distribution of $X_\xi$ is a generalized Gamma convolution (GGC) and is hence SD - see Example 4.3.4 in \cite{B}. The natural question whether $X_\xi$ is SD or even ID for $\xi \in (-1,0)$ was first raised by Steutel in 1973 - see Section 3.4 in \cite{St}, and has remained open ever since. In section 4.5 of \cite{B} - see also the Appendix B.3 of \cite{SVH}, this problem is rephrased in the broader context  of generalized Gamma distributions. The latter are power transformations of $\Ga_t$ where $\{\Ga_t, \, t\ge 0\}$ is the Gamma subordinator, and  can be thought of as size-biased samplings of $X_\xi$ when $\xi <0$, in view of the formul\ae
$$\EE[f(\Ga_t^\xi)]\; =\; \frac{\EE[f(X_\xi)X_\xi^{u_{}}]}{\EE[X_\xi^{u_{}}]}$$
valid for every $f$  bounded continuous and $t > 0,$ with $u = (t-1)/\xi.$ Recall in passing that Steutel's equation - see e.g. Theorem 51.1 in \cite{S1} - establishes a precise link between size-biased sampling of order one  and infinite divisibility for integrable positive random variables. In this note, we provide an answer to the above questions of \cite{St, B}.

\begin{THE} For every $\xi \in (-1,0) $ and $t > 0,$ the random variable $\Ga_t^\xi$ is {\rm SD}.
\end{THE}

As a direct consequence of this result, all Fr\'echet distributions are SD and the extreme value distribution $F_\xi$ is ID if and only if $\xi \not\in (0,1).$ Contrary to the case $\vert \xi\vert \ge 1,$ our argument is probabilistic and consists in showing that $\Ga_t^\xi$ is distributed as the exponential functional of some spectrally negative L\'evy process. This extends a classical result of Dufresne \cite{D} for the case $\xi = -1.$ The identification is made possible thanks to a entire moment method due to Bertoin-Yor \cite{BY2}, which applies in our context as a case study. The proof is given in the next section. 

In Section 3, we review the possible interpretations of the infinite divisibility of $\Ga_t^\xi$ for $\xi < 0.$ The classical case $\xi = -1$ allows at least four different formulations in terms of processes, and also an explicit computation of the L\'evy density which shows the GGC property without the HCM argument. For $\xi <- 1$ the ID property is only known by analytical means and there is no direct probabilistic explanation, save for the case $t = 1$ by subordination or, tentatively, the spectral theory of a certain spectrally positive Markov processes. The situation for $\xi \in (-1,0)$ is exactly the opposite since in addition to the exponential functional argument, the ID property can also be obtained rigorously by a first-passage time argument for a spectrally positive Markov processes. On the other hand there is no analytic proof of the ID property for $\xi \in(-1,0)$. In this situation the GGC character of $\Ga_t^\xi$ remains in particular an open question, which we plan to tackle in some further research.

\section{Proof of the Theorem}

We begin with a computation on the Gamma function.

\begin{LEM} For every $\a\in(0,1)$ and $u, t > 0$ one has
$$\frac{u\Ga(t + \a (u+1))}{\Ga(t+\a u)}\; =\; \lpa\frac{\Ga(t + \a)}{\Ga(t)}\rpa u \; +\; \int_{-\infty}^0 (e^{u x} - 1 - ux)f_{\a,t}(x)dx,$$
where
$$f_{\a,t}(x)\; =\; \frac{e^{(1+t/\a)x}(\a+e^{x/\a} +t(1-e^{x/\a}))}{\a\Ga(1-\a)(1-e^{x/\a})^{\a+2}}$$
is the density of a L\'evy measure on $(-\infty, 0).$
\end{LEM}

\proof We set $\lbd = t +\a u > 0$ and compute
\begin{eqnarray*}
\frac{\Gamma(\lbd +\a)}{\Gamma(\lbd)} & = & \frac{\lbd \beta (\lbd +\a, 1-\a)}{\Gamma(1-\a)}\\
& =&  \frac{\lbd}{\Gamma(1-\a)}\int_0^{+\infty} \!\!\frac{e^{-(\a +\lbd) x}}{(1-e^{-x})^{\a}}\, dx \\
& =&  \frac{\a}{\Gamma(1-\a)}\int_0^{+\infty} (1- e^{-\lbd x}) \frac{e^{-\a x}}{(1-e^{-x})^{\a +1}}\, dx \\
\end{eqnarray*}
where the second equality comes from a change of variable and the third from an integration by parts. This yields
\begin{eqnarray*}
\frac{u\Ga(t + \a (u+1))}{\Ga(t+\a u)} & =&  \frac{\a u}{\Gamma(1-\a)}\int_0^{+\infty} (1- e^{-(t+\a u) x}) \frac{e^{-\a x}}{(1-e^{-x})^{\a +1}}\, dx \\
& = & \lpa\frac{\Ga(t + \a)}{\Ga(t)}\rpa u \; +\; \frac{\a u}{\Gamma(1-\a)}\int^0_{-\infty} (1- e^{\a u x}) \frac{e^{(\a +t) x}}{(1-e^x)^{\a +1}}\, dx\\
& = & \lpa\frac{\Ga(t + \a)}{\Ga(t)}\rpa u \; +\; \int_{-\infty}^0 (e^{u x} - 1 - ux)f_{\a,t}(x)dx
\end{eqnarray*}
where again, the second equality comes from a change of variable and the third from an integration by parts.

\endproof

\begin{REMS} {\em (a) The above proof follows \cite{BY1} p. 102. Notice in passing that some computations performed in \cite{BY1} are slightly erroneous. For example the subordinator whose exponential functional is distributed as $\tau_\a^{-\a}$ (with the notation of \cite{BY1}) has no drift, but it is also killed at rate $1/\Gamma(1-\a).$\\

\noindent
(b) The above decomposition extends to $\a = 1$ since
$$\frac{u\Ga(t + (u+1))}{\Ga(t+ u)}\; =\; u(t+u)$$
is the L\'evy-Khintchine exponent of a drifted Brownian motion (the latter was already noticed in \cite{BY2} - see Example 3 therein - in order to recover Dufresne's identity). However, such a formula does not seem to exist for $\a > 1.$}
\end{REMS}

\bigskip

\noindent
{\bf End of the proof.} Fix $\xi \in(-1,0), t > 0,$ and set $\a = -\xi\in(0,1)$ for simplicity. The entire moments of $\Ga_t^{\a}$ are given for every $n\ge 1$ by 
\begin{eqnarray*}
\EE[\Ga_t^{\a n}] & = & \frac{\Ga (t+\a n)}{\Ga (t)} \\
& = & \frac{\Ga (t+\a)}{\Ga (t)}\times\cdots\times \frac{\Ga (t+\a n)}{\Ga (t+\a (n-1))}\; =\; m\, \frac{\psi(1)\ldots\psi(n-1)}{(n-1)!}
\end{eqnarray*}
with the notation
$$\psi(u)\; =\; \frac{u\Ga(t + \a (u+1))}{\Ga(t+\a u)}\; =\; \lpa\frac{\Ga(t + \a)}{\Ga(t)}\rpa u \; +\;\int_{-\infty}^0 (e^{u x} - 1 - ux)f_{\a,t}(x)dx$$
by the Lemma, and
$$m \; = \; \frac{\Ga(t + \a)}{\Ga(t)}=\psi'(0+).$$
It is clear that $\psi$ is the L\'evy-Khintchine exponent of a spectrally negative L\'evy process $Z$ with infinite variation and mean $m >0.$ By Proposition 2 in \cite{BY2}, this entails 
$$\EE[\Ga_t^{\a n}] \; = \; \EE[I^{-n}]$$
for every $n \ge 1,$ where $I$ is the exponential functional of $Z$:
$$I\; =\; \int_0^\infty e^{-Z_s} \, ds.$$ 
Since $Z$ has no positive jumps, Proposition 2 in \cite{BY2} shows also that the random variable $1/I$ is moment-determinate, whence
$$\Ga_t^\xi\;\elaw\; I.$$
The self-decomposability of $I$ is a direct consequence of the Markov property. More precisely, introducing the stopping-time $T_y = \inf\{s > 0,\; Z_s = y\}$ for every $ y >0$, the fact that $Z_s \to +\infty$ a.s. as $s\to +\infty$ and the absence of positive jumps entail that $T_y < +\infty$ a.s. Decomposing, we get
$$I \; =\;  \int_0^{T_y} e^{-Z_s} \, ds\; +\;  \int_{T_y}^\infty e^{-Z_s} \, ds\; \elaw\;  \int_0^{T_y} e^{-Z_s} \, ds\; +\; e^{-y} \int_0^\infty e^{-Z'_s} \, ds$$
where $Z'$ is an independent copy of $Z$ and the second equality follows from the Markov property at $T_y$. This shows that for every $c\in (0,1)$ there is an independent factorization
$$I\; =\; cI\; +\; I_c$$ 
for some random variable $I_c,$ in other words that  $I\elaw \Gamma_t^\xi$ is self-decomposable.

\qed

\begin{REMS} {\em (a) By the above Remark 1 (b), it does not seem that $\Ga_t^\xi$ is distributed as the exponential functional of a L\'evy process for $\xi <- 1.$ It would be interesting to have an explanation of the infinite divisibility of $\Ga_t^\xi$ in terms of processes when $\xi<- 1.$ See next section for a more precise conjecture in the case $t=1$. \\

\noindent
(b) The self-decomposability of $S_\a^{\a_{}}$ for every $\a\in (0,1)$  has been shown by Patie \cite{Patou} in using the same kind of argument. Specifically, one can write
$$\EE[\Sa^{n\a}]\; =\; \frac{\Gamma(1+n)}{\Gamma(1+\a n)}\; =\; m\, \frac{\psi(1)\ldots\psi(n-1)}{(n-1)!}$$
where we use the same notation as above and, correcting small mistakes made in Paragraph 3.2 of \cite{Patou},
$$\psi(u)\; =\; \frac{u}{\Ga(1+\a)} \; +\; \int_{-\infty}^0 (e^{u x} - 1 - ux) \frac{(1-\a)e^{x/\a}((2-\a)e^{x/\a} +(1-e^{x/\a}))}{\a^2\Ga(1+\a)(1-e^{x/\a})^{3-\a}}\, dx$$
is the L\'evy-Khintchine exponent of some spectrally negative L\'evy process with positive mean. Setting $\a = t = 1/2$ and comparing the above expression to the one in the Lemma, one can check the well-known identity in law
\begin{equation}
\label{TheID}
\sqrt{S_{1/2}}\; =\; \frac{1}{2\sqrt{\Ga_{1/2}}}\cdot
\end{equation}
The present paper shows that all positive powers of $S_{1/2}$ are SD and one may wonder if the same is true for $\Sa$ with any $\a\in (0,1).$ See \cite{JS} for related results and also for a characterization of the SD property of negative powers of $\Sa$ when $\a\leq 1/2.$}
\end{REMS}

\section{Further remarks and open questions}

In this section we would like to review several existing or tentative approaches for the ID, SD and GGC properties of the distribution of $\Ga_t^{\xi}$ or $X_\xi$ when $\xi\le 0.$ 

\subsection{The case $\xi = 0$} This is a rather specific situation but we include it here for completeness. As mentioned in Section 3.4 of \cite{St}, the SD property of the two-sided $X_0$ is a direct consequence of the extreme value theory because
\begin{equation}
\label{EVT}
L_1+\frac{L_2}{2} +\cdots +\frac{L_n}{n} -\log n\; \elaw\;\max(L_1, \ldots, L_n) - \log n\; \claw\; X_0
\end{equation}
as $n\to +\infty,$ where $L_1,\ldots, L_n$ are independent copies of $L\sim$ Exp $(1).$ The above identity and convergence in law, readily obtained in comparing Laplace transforms and distribution functions, yield after some further computations the following closed expression for the Laplace transform of $X_0$:
$$\EE[e^{-\lambda X_0}]\; =\; \Gamma(1+\lambda)\; =\; \exp\lcr -\gamma\lambda \, +\, \int_0^\infty (e^{-\lambda x}-1+\lambda x)\frac{dx}{x(e^x-1)}\rcr,$$
where $\gamma$ is Euler's constant. The complete monotonicity of $1/(e^x -1)$ shows then that $X_0$ is an extended GGC in the sense of Chapter 7 in \cite{B}. See also Exercise 18.19 in \cite{S1} and Example 7.2.3 in \cite{B} for another argument based on Pick functions, recovering (\ref{EVT}).

\subsection{The case $\xi = -1$} This is the classical situation, very well-known, but we give some details for comparison purposes. The ID property of $X_{-1}$ can first be understood by the sole fact that
\begin{equation*}
\lim_{n\to +\infty} \lpa \frac{nx}{1+nx}\rpa^n\; =\; e^{-1/x}
\end{equation*}
because the left-hand side is the first-passage time distribution function of a certain birth and death process - see Theorem 3.1 and (3.3) in \cite{St}. The random variable $\Ga_t^{-1}$ is also a GIG and is hence distributed as the unilateral first-passage time of a diffusion \cite{BBH}, which explains its  infinite divisibility by continuity and the Markov property. More precisely one has
\begin{equation}
\label{Bess}
\frac{1}{4\Ga_t}\; \elaw\; \inf\{u > 0,\; X^t_u\, =\, 0\}
\end{equation}
where $\{X^t_u, \, u\ge 0\}$ is a Bessel process of dimension $2(1-t)$ starting from one. The SD property follows as above from  Dufresne's identity \cite{D}, which reads
\begin{equation}
\label{Duf}
\frac{2}{\Ga_t}\; \elaw\; \int_0^\infty e^{B_u -tu/2} \, du
\end{equation}
where $\{B_u, \, u\ge 0\}$ is a standard linear Brownian motion. Also,  Exercise 16.4 in \cite{S1} shows that $\Ga_t^{-1}$ is the one-dimensional marginal of a certain self-similar additive process, whence its self-decomposability by Theorem 16.1 in \cite{S1}. The link between this latter interpretation and (\ref{Bess}) and (\ref{Duf}) has been explained in depth in \cite{Y}. 

It does not seem that these four interpretations can provide any explicit  information on the L\'evy-Khintchine exponent of $\Ga_t^{-1}.$ But in this case analytical computations can also be performed. More precisely, taking for simplicity the same normalization as in (\ref{Bess}) and setting $\varphi_t(\lbd) = -\log \EE[e^{-\lbd/4\Ga_t}],$  formul\ae\, (7.12.23), (7.11.25) and (7.11.26) in \cite{E} entail
\begin{equation}
\label{phip}
\varphi_t'(\lbd)\; =\; \frac{K_{t-1}(\sqrt{\lbd})}{2\sqrt{\lbd} K_t(\sqrt{\lbd})}
\end{equation}
where $K_t$ is the Macdonald function. This shows $\varphi_{1/2}'(\lbd) = 1/2\sqrt{\lbd}$ viz. $\varphi_{1/2}(\lbd) = \sqrt{\lbd}$ when $t=1/2,$ and one recovers the identity (\ref{TheID}). For $t=3/2,$ one obtains
$$2\varphi_{3/2}'(\lbd) \; =\; \frac{1}{1+\sqrt{\lbd}}\; =\; \EE[e^{-\lbd (L^2\times S_{1/2})}]\; =\; \int_0^\infty \lpa \frac{1}{\lbd +x}\rpa \frac{\sqrt{x}\, dx}{\pi (1+x)}$$
where the first equality follows from Formula (7.2.40) in \cite{E}, and the third equality from Exercise 29.16 in \cite{S1} and (2.2.5) in \cite{B}.
This means precisely - see (3.1.1) in \cite{B} - that the distribution of $1/4\Gamma_{3/2}$ is a GGC with zero drift and Thorin measure 
$$U_{3/2}(dx)\; =\; \frac{\sqrt{x}\, dx}{2\pi (1+x)}\cdot$$
The latter property can be extended to {\em all} values of $t$ thanks to a result originally due to Grosswald \cite{Gr} on Student's distribution. Together with (\ref{phip}), the main theorem in \cite{Gr} entails namely that the distribution of $1/4\Gamma_t$ is a GGC with zero drift and Thorin measure 
$$U_t(dx)\; =\; \frac{1}{\pi^2x (J^2_t(\sqrt{x}) + Y^2_t (\sqrt{x}))}$$
where $J_t$ and $Y_t$ are the usual Bessel functions of the first kind - see \cite{E} p. 4.

\subsection{The case $\xi\in (-1,0)$} In this situation, the present paper yields an interpretation of the self-decomposability of $\Ga_t^\xi$ by the identification
$$\Ga^{\xi}_t\; \elaw\; \int_0^\infty e^{-Z_u} \, du,$$
where $Z$ is a spectrally negative L\'evy process. Another explanation,  similar to (\ref{Bess}), can then be obtained by the Lamperti transformation - see e.g. the introduction of \cite{BY2} for an account and references. More precisely, setting
$$Y_u\; =\; \exp[ -Z_{\tau_u}]$$
with the notation $\tau_u = \inf\{ s > 0, \; \int_0^s e^{-Z_v} \, dv \, > u\}$, then $Y = \{ Y_u, \; 0\le u <\Ga^{\xi}_t\}$ is a spectrally positive Markov process (which is also self-similar) starting from one and we have
$$\Ga^{\xi}_t\; \elaw\;  \inf\{ u > 0, \; Y_u = 0\},$$
so that the infinite divisibility of $\Ga^{\xi}_t$ (but not, or at least not directly, its self-decomposability) is a consequence of the Markov property and the absence of negative jumps for $Y.$ It would be interesting to see if $\Ga^{\xi}_t$ could be embedded in some self-similar additive process analogous to the Brownian escape process of the case $\xi =1,$ described in Exercise 16.4 of \cite{S1}. 

Our main result can also be interpreted analytically in terms of generalized Bessel functions. Setting $\a = -\xi$ and writing down
\begin{equation}
\label{GB}
\EE[e^{-\lbd \Ga^{\xi}_t}]\; =\; \frac{1}{\a\Ga(t)}\int_0^\infty x^{-t\a -1} e^{-\lbd x +x^{-1/\a}} dx\; =\; \frac{ Z_{1/\a}^{t/\a}(\lbd)}{\a\Ga(t)}
\end{equation}
with the notation (1.7.42) of \cite{Ki}, the infinite divisibility of $\Ga^{\xi}_t$ entails that the function
\begin{equation}
\label{GBF}
\lbd\; \mapsto\; - \lpa\frac{Z_{\rho}^{\nu\,'}(\lbd)}{Z_\rho^\nu(\lbd)}\rpa
\end{equation}
is completely monotone for any $\rho > 1$ and $\nu > 0.$ One might ask if the latter function is also a Stieltjes transform, which is equivalent to the GGC property for the distribution of $\Ga^{\xi}_t$ - see Chapter 3 in \cite{B}. Indeed, it is very natural to conjecture such a property for $\xi\in (-1,0)$ in view of the above cases $\xi = 0$ and $\xi=-1.$ Compared to classical Bessel functions, the theory of generalized Bessel functions is however rather incomplete, and proving like in \cite{Gr} that the function (\ref{GBF}) is a Stieltjes transform is believed to be challenging.
 
\subsection{The case $\xi <- 1$} In this situation the GGC property of the distribution of $\Ga_t^\xi$ is most quickly obtained from the HCM character of the density function - see Chapter 5 and especially Example 5.5.2 in \cite{B}. Notice that this analytical argument extends to $\xi = -1$ but not to $\xi\in(-1,0)$ since otherwise $\Ga_t^{-\xi}$ would also have a HCM density and hence be ID, which is false. This entails that the function in (\ref{GBF}) is indeed a Stieltjes transform for any $\rho \in (0,1)$ and $\nu > 0,$ and it would be  interesting to identify the underlying Thorin measure as in Grosswald's theorem.

A probabilistic interpretation of the self-decomposability of $\Ga_1^\xi = L^\xi$ can also be given by Bochner's subordination. Setting $\a = -1/\xi\in (0,1),$ one has indeed 
$$L^\xi\;\elaw\; L^{-1}\, \times\, \Sa\; \elaw\; S^\a_{L^{-\a}}.$$
where $\{S^\a_u, \, u\ge 0\}$ stands for the $\a-$stable subordinator with marginal $S^\a_1 \elaw \Sa.$ Since $L^{-\a}$ is SD by our result, this means that $L^\xi$ is the marginal of some subordinator which is itself subordinated to the $\a-$stable one, and Proposition 4.1. in \cite{S2} shows that $L^\xi$ is SD. Besides, setting $\varphi_\a$ resp. $\varphi_\xi$ for the L\'evy-Khintchine exponent of $L^{-\a}$ resp. $L^\xi,$ one deduces from Theorem 30.4 in \cite{S1} the following relationship
$$\varphi_\xi(\lbd)\; =\; \varphi_\a(\lbd^\a).$$
Another, tentative, probabilistic interpretation of the self-decomposability of $L^\xi$ could be given in terms of a certain spectrally positive Markov process. Setting $\a = -1/\xi$ and $y_\a(\lbd) = \EE[e^{-\lbd L^\xi}],$ Theorem 4.17 p. 258 in \cite{Ki} and (\ref{GB}) above show that $y_\a$ 
is a solution to the fractional differential equation
$$x D_-^{\a+1} y_\a \; -\; \a y_\a \; =\; 0,$$
where $D_-^{\a+1}$ is a fractional Riemann-Liouville derivative - see Section 2.1 in \cite{Ki}. When $\a = 1$ viz. $\xi = -1$ the above amounts to a Bessel equation and Feller's theory applies, making $L^{-1}$ the first-passage time of a Bessel process of index 0 - see \cite{Ke}. When $\a \in (0,1)$ the operator $D_-^{\a+1}$ is the infinitesimal generator of a spectrally positive $(1+\a)-$stable L\'evy process reflected at its minimum, which is a spectrally positive Markov process - see Section 3 in \cite{PS} and the references therein. By downward continuity, one may wonder if $L^\xi$ cannot be viewed as the first-passage time of some scale-transformation of the latter, eventhough no Feller's theory is available for fractional operators whose order lies in $(1,2).$

The above probabilistic interpretations do not seem to hold for $t\neq 1.$ On the one hand, Theorem 4.17 in \cite{Ki} yields then an equation with {\em two}  fractional derivatives of different order for the Laplace transform of $\Ga_t^\xi$. On the other hand, keeping the notation $\a = -1/\xi,$ Theorem 1 in \cite{SHS} shows the factorization 
$$\Ga_t^\xi\; =\; \Ga_{\a t}^{-1}\,\times\, \Sa^{(t)}$$
where $\Sa^{(t)}$ is the size-biased sampling of $\Sa$ of order $-\a t,$ viz.
$$\EE[f(\Sa^{(t)})]\; =\;  \frac{\EE[f(\Sa)\Sa^{-\a t}]}{\EE[\Sa^{-\a t}]}\cdot$$
The GGC character of $\Sa^{(t)}$ follows from that of $\Sa$ and Theorem 6.2.4. in \cite{B}, which shows by the above case $\xi=-1$ that 
$\Ga_t^\xi$ is the independent product of two SD random variables. By Theorem 16.1 in \cite{S1}, this entails that there exist two independent 1-self-similar additive increasing processes $Y$ and $Z$ such that
$$\Ga_t^\xi\; \elaw \; Y_{Z_1}.$$
Unfortunately, contrary to Bochner's subordination the independent composition of two additive processes is not necessarily an additive process anymore, so that the above identity does not provide a probabilistic proof of the self-decomposability of $\Ga_t^\xi.$ 

\bigskip

\noindent
{\bf Acknowledgements.} The authors are grateful to L.~Bondesson for the interest he took in this work, and several useful comments. Ce travail a b\'en\'efici\'e d'une aide de l'Agence Nationale de la Recherche portant la r\'ef\'erence ANR-09-BLAN-0084-01.

\end{document}